\documentclass[11pt]{article}      % 
\usepackage{amsmath,amssymb, amsthm, eucal,bm,accents,mathtools}  
\usepackage{graphicx, amscd}
\usepackage{framed}

\usepackage{newtxtext,newtxmath} %<<<<< seems better than mathptmx
\usepackage[T1]{fontenc}     %for better hyphenation etc

\usepackage{scrextend}     %used to change margin for a part of doc 
\usepackage{enumerate}       %used for better lists
\usepackage{float}                 % use in \begin{figure}[H] to enforce position
\usepackage{tikz}
\usetikzlibrary{matrix, arrows, calc}
\usetikzlibrary{arrows,shapes,positioning}

\theoremstyle{plain}
\newtheorem{theorem}{Theorem}[section]            %{commend}{in print}[numbered within]       
\newtheorem{proposition}[theorem]{Proposition}  %{commend}[numbered with]{in print as}

\theoremstyle{definition}

\numberwithin{theorem}{section}
\numberwithin{equation}{section}
\numberwithin{figure}{section}
\newcommand{\gaction}[2]{\genfrac{}{}{0.5pt}{}{#1}{#2}% 
                        \!\lower2pt\hbox{\rotatebox[origin=c]{-90}{{$\looparrowright$}}}}
\newcommand{\dotfraction}[2]{\genfrac{}{}{0.5pt}{}{#1}{#2}% 
                        \!\lower.5pt\hbox{{$\circ$}}}

\usepackage{titlesec}                                                           %fixes section titles
\titlelabel{\thetitle.\ }                                                         %number
\titleformat*{\section}{\fontsize{14pt}{14pt} \bf}        %size
\usepackage{fancyhdr}
\usepackage{sidecap}  %for side caption of figures.  Use :
\usepackage[hang,small,sc]{caption}

\def\-{\hbox{\raisebox{.75pt}{-}}}

\usepackage{fancybox}

\def\smalll{\scriptsize}
\begin{document}

%\title{Area of coronas in the Apollonian disk packing}
\title{\huge Fibonacci numbers and Ford circles}

\author{Jerzy Kocik %\thanks{support}\\
    \\ \small Department of Mathematics, Southern Illinois University, Carbondale, IL62901
   \\ \small jkocik{@}siu.edu
}

%\date{\small (\the\day \  \the\month \ \the\year)}
%\date{}  % NO DATE

\date{\small\today}

\maketitle

\begin{abstract}
\noindent
An amusing connection between Ford circles, Fibonacci numbers, and golden ratio is shown.
Namely, 
certain tangency points of Ford circles are concyclic and involve Fibonacci numbers.
They form four circles that cut the x-axis at points related to the golden ratio.
%four different circles of passing through point of tangecy of certain Fod circles circles 
\\[5pt]
{\bf Keywords:} 
Fibonacci numbers, Ford circles, golden ratio.
%?Inspection, exploratory inspection
\\[5pt]
\scriptsize {\bf MSC:} 11B39, %  	Fibonacci and Lucas numbers and polynomials and generalizations
                                %51M04, %Elementary problems in E Geo
                                52C26,  %Circle packings and discrete conformal geometry
                                28A80,  %Fractals
                                51M15. %geometric constructions
                                   %11A99  % Number theory: None of the above, but in this section

\end{abstract}

%cute
%pleasing
%delectable

%--------------------------------------------------------------------------------------
\section{Introduction}

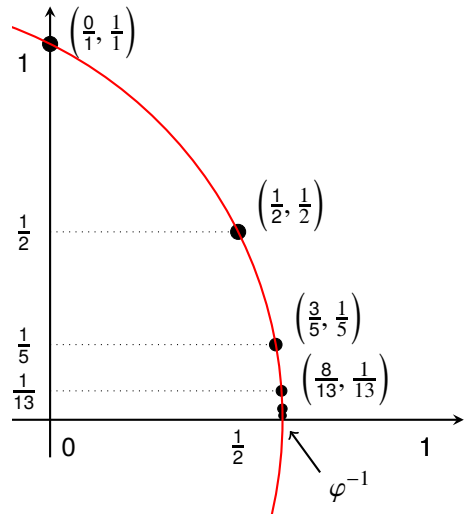
\begin{figure}[H]
\centering
\begin{tikzpicture}[scale=5]
\clip (-.1, -.25) rectangle (1.1, 1.2);
%%-----------------------  
\draw [-stealth, thick] (-1,0)--(1.1,0);  %line
\draw [-stealth, thick] (0,-.1)--(0,1.1);  %line

\draw [dotted] (0,1/2)--(1/2,1/2);  %line
\draw [dotted] (0,1/5)--(3/5,1/5);  %line
\draw [dotted] (0,1/13)--(8/13,1/13);  %line

%\foreach \a /  \b  in {
%0/1, 1/1,
%1/2,
%1/3,2/3,
%1/4,3/4,
%1/5,2/5,3/5,4/5,
%1/6,1/6,
%1/7,2/7,3/7,4/7,5/7,6/7,
%1/8,3/8,5/8,7/8
%}
%\draw [thick]  (\a/\b, 1/2*1/\b*1/\b) circle (1/2*1/\b*1/\b)
%%   (1+\a/\b, 1/2*1/\b*1/\b) circle (1/2*1/\b*1/\b)
%   (-1+\a/\b, 1/2*1/\b*1/\b) circle (1/2*1/\b*1/\b)
%   (-2+\a/\b, 1/2*1/\b*1/\b) circle (1/2*1/\b*1/\b);

%\draw  [->,  red, bend right=-20] (1/2,1/2) to  (3/5, 1/5);

\draw [fill=black] (1/2,1/2) circle (2/100); %diagonal
\draw [fill=black] (3/5,1/5) circle (2/120); %diagonal
\draw [fill=black] (8/13,1/13) circle (2/140); %diagonal
\draw [fill=black] (21/34,1/34) circle (2/160); %diagonal
\draw [fill=black] (55/89,1/89) circle (2/200); %diagonal

\draw [fill=black] (0,1) circle (2/100); %diagonal
%\draw [fill=black] (-1,1) circle (2/120); %diagonal
%
%\draw [fill=black] (-1-1/2,1/2) circle (2/100); %diagonal
%\draw [fill=black] (-1-3/5,1/5) circle (2/120); %diagonal
%\draw [fill=black] (-1-8/13,1/13) circle (2/140); %diagonal
%\draw [fill=black] (-1-21/34,1/34) circle (2/160); %diagonal
%\draw [fill=black] (-1-55/89,1/89) circle (2/200); %diagonal

%\draw [fill=black] (1/2,1/2) circle (1/100); %diagonal
%\draw [fill=black] (2/5,1/5) circle (1/120); %diagonal
%\draw [fill=black] (5/13,1/13) circle (1/140); %diagonal
%\draw [fill=black] (13/34,1/34) circle (1/160); %diagonal
%\draw [fill=black] (34/89,1/89) circle (1/200); %diagonal

\node at (.05, -.07) [scale=1, color=black] {$\sf 0$};
\node at (1, -.07) [scale=1, color=black] {$\sf 1$};
\node at (0-.07, .95) [scale=1, color=black] {$\sf 1$};
\node at (1/2, -.07) [scale=1, color=black] {$\frac{\sf 1}{\sf 2}$};
\node at (1, -.07) [scale=1, color=black] {$\sf 1$};
\node at (0-.07, 1/2) [scale=1, color=black] {$\frac{\sf 1}{\sf 2}$};
\node at (0-.07, 1/5) [scale=1, color=black] {$\frac{\sf 1}{\sf 5}$};
\node at (0-.07, 1/13) [scale=1, color=black] {$\frac{\sf 1}{\sf 13}$};

\node at (0/1+.14, 1/1+.03) [scale=1, color=black] {$\left(\frac{\sf 0}{\sf 1}, \frac{1}{1}\right)$};
\node at (1/2+.14, 1/2+.07) [scale=1, color=black] {$\left(\frac{\sf 1}{\sf 2}, \frac{1}{2}\right)$};
\node at (3/5+.14, 1/5+.07) [scale=1, color=black] {$\left(\frac{\sf 3}{\sf 5}, \frac{1}{5}\right)$};
\node at (8/13+.17, 1/13+.03) [scale=1, color=black] {$\left(\frac{\sf 8}{\sf 13}, \frac{1}{13}\right)$};
%\node at (2/3, -.07) [scale=1, color=black] {$\frac{\sf 2}{\sf 3}$};
%\node at (1/4, -.07) [scale=1, color=black] {$\frac{\sf 1}{\sf 4}$};
%\node at (3/4, -.07) [scale=1, color=black] {$\frac{\sf 3}{\sf 4}$};
%\node at (1/5, -.07) [scale=1, color=black] {$\frac{\sf 1}{\sf 5}$};
%\node at (2/5, -.07) [scale=1, color=black] {$\frac{\sf 2}{\sf 5}$};
%\node at (3/5, -.07) [scale=1, color=black] {$\frac{\sf 3}{\sf 5}$};
%\node at (4/5, -.07) [scale=1, color=black] {$\frac{\sf 4}{\sf 5}$};
%\node at (1/6, -.07) [scale=1, color=black] {$\frac{\sf 1}{\sf 6}$};
%
%\node at (-1/1, -.07) [scale=1, color=black] {-$\frac{\sf 1}{\sf 1}$};
%\node at (-1/2, -.07) [scale=1, color=black] {-$\frac{\sf 1}{\sf 2}$};

%\draw [->,thick] (.382, -.14) -- (.382, -.01);
\draw [->,thick] (.618+.1, -.14) -- (.618+.02, -.03); 

%\node at (1/2.618, -.177) [scale=1, color=black] {$\ \varphi^{-2}$};
\node at (1/1.618+.17, -.177) [scale=1, color=black] {$\ \varphi^{-1}$};

\draw [thick, red] (-1/2,0) circle (1/2+.618); %circle

%\draw [thick, blue] (3/2,0) circle (1/2+.618); %circle

\end{tikzpicture}
\caption{Fibonacci points}
\label{fig:points}
\end{figure}

Let $F_n$ be the Fibonacci sequence labeled as follows:
$$
F_0=0, \quad F_1=1,\quad F_{n+1}= F_n+F_{n-1}, \quad n\in \mathbb Z
$$
It looks and feels ``1-dimensional''.
But there is a hidden circle in this sequence.
Consider ratios of consecutive pairs, and points with the following coordinates:
%$$
\begin{equation}
\label{eq:points}
\left( \frac{F_{2n}}{F_{2n+1}},\,\frac{1}{F_{2n+1}}\right)
\end{equation}
%$$
The $x$-coordinates of \eqref{eq:points}  are ratios produced from Fibonacci numbers in the following manner:
$$%\hspace{-.5in}
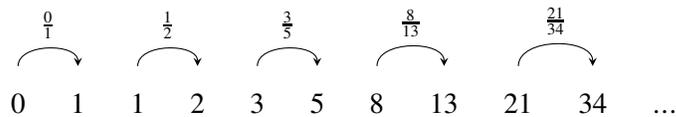
\begin{tikzpicture}[baseline=-0.7ex]
    \matrix (m) [ matrix of math nodes,
                         row sep=1.7em,
                         column sep=.9em,
                         text height=2.8ex, text depth=1.1ex] 
      {
                    0   
                &  1  
                &  1  
                &  2  
                &  3                  
                &  5                 
                &  8             
                &  13 
                &  21                  
                &  34
                & ...
                \\
     };

    \path[-stealth]
%horizontal
        (m-1-1) edge [out=90, in=90]  node[above ] { \smalll $\frac{0}{1}$}  (m-1-2)
        (m-1-3) edge [out=90, in=90]  node[above ] { \smalll $\frac{1}{2}$}  (m-1-4)        
        (m-1-5) edge [out=90, in=90]  node[above ] { \smalll $\frac{3}{5}$}  (m-1-6)        
        (m-1-7) edge [out=90, in=90]  node[above ] { \smalll $\frac{8}{13}$}  (m-1-8)
         (m-1-9) edge [out=90, in=90]  node[above ] { \smalll $\frac{21}{34}$}  (m-1-10)
    
;
\end{tikzpicture}   
%\end{equation}
$$
It turns out that points \eqref{eq:points} are concyclic!  
See Figure \ref{fig:points}.
But there is more.  
The so-defined  circle passes through the points of tangency of Ford circles.  
Actually, there are four such +circles of similar connection to Fibonacci numbers.

%--------------------------------------------------------------------------------
\section{Ford circles}

Let us recall some basic facts about Ford circles.

\begin{proposition}
\label{thm:Ford}{\rm [L. Ford]}
\sf
Define
$
K[p,m]%$ = \hbox{circle of radius} \frac{1}{2m^2} \hbox{at point } x=\frac{p}{m}
$
to be a circle of radius $1/2m^2$ tangent to x-axis at point $x=\frac{p}{m}$.
Then 
%$$
\begin{equation}
\label{eq:matrix}
\det\begin{bmatrix}p&q\\m&n\end{bmatrix} = \pm1
\qquad \Rightarrow\qquad
K[p,m] \ \ \hbox{and} \ \   K[q,n] \ \hbox {are tangent}
\end{equation}
%$$
and otherwise are disjoint.
Moreover, the Ford circle at 
%$$ 
\begin{equation}
\label{eq:Farey}
x = \frac{p}{n} \oplus \frac {q}{m} \equiv \frac{p+q}{m+n}
\end{equation}
%$$ 
is inscribed in the triangle-like region between the two tangent circles 
at $p/n$ and $q/m$, and the $x$-axis.
\end{proposition}

\noindent
For the proof see \cite{Ford}.

~\\
{\bf Remark:} We shall also include the horizontal line $y=1$ among the Ford circles, 
as a line over the point at infinity,  $x=\frac{1}{0}$.

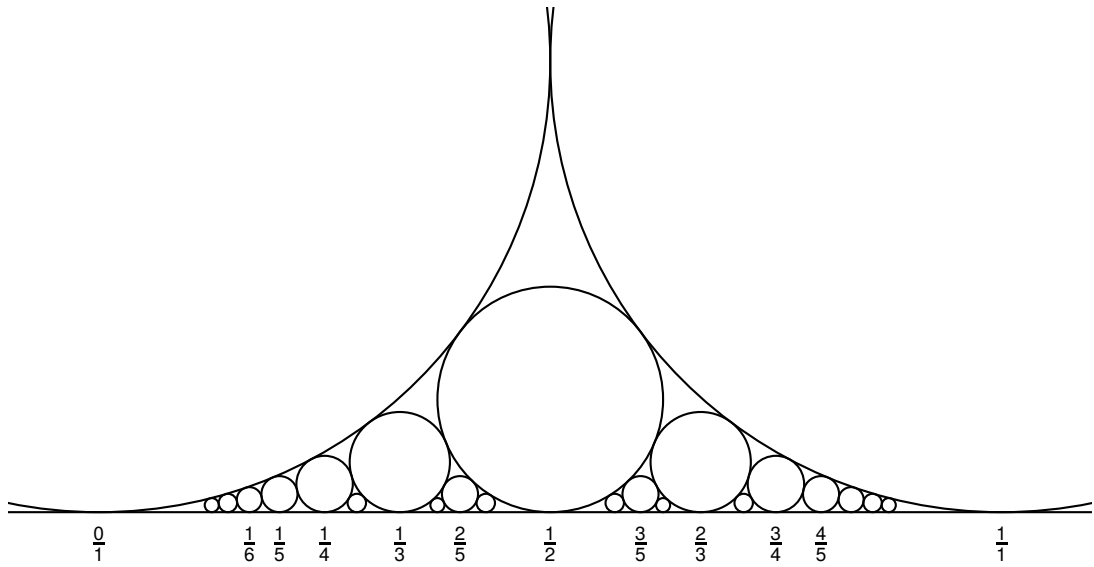
\begin{figure}[H]
\centering
\begin{tikzpicture}[scale=12]
\clip (-.1, -.12) rectangle (1.1, .56);
%%-----------------------  
\draw [thick] (-.2,0)--(1.2,0);  %line

\foreach \a /  \b  in {
0/1, 1/1, 1/2, 1/3,2/3, 1/4,3/4, 1/5,2/5,3/5,4/5,
1/6,5/6, 1/7,2/7,3/7,4/7,5/7,6/7, 1/8,3/8,5/8,7/8
}
\draw [thick]  (\a/\b, 1/2*1/\b*1/\b) circle (1/2*1/\b*1/\b)
%   (1+\a/\b, 1/2*1/\b*1/\b) circle (1/2*1/\b*1/\b)
%   (-1+\a/\b, 1/2*1/\b*1/\b) circle (1/2*1/\b*1/\b)
%  (-2+\a/\b, 1/2*1/\b*1/\b) circle (1/2*1/\b*1/\b)
;

\node at (0/1, -.035) [scale=1, color=black] {$\frac{\sf 0}{\sf 1}$};
\node at (1/1, -.035) [scale=1, color=black] {$\frac{\sf 1}{\sf 1}$};
\node at (1/2, -.035) [scale=1, color=black] {$\frac{\sf 1}{\sf 2}$};
\node at (1/3, -.035) [scale=1, color=black] {$\frac{\sf 1}{\sf 3}$};
\node at (2/3, -.035) [scale=1, color=black] {$\frac{\sf 2}{\sf 3}$};
\node at (1/4, -.035) [scale=1, color=black] {$\frac{\sf 1}{\sf 4}$};
\node at (3/4, -.035) [scale=1, color=black] {$\frac{\sf 3}{\sf 4}$};
\node at (1/5, -.035) [scale=1, color=black] {$\frac{\sf 1}{\sf 5}$};
\node at (2/5, -.035) [scale=1, color=black] {$\frac{\sf 2}{\sf 5}$};
\node at (3/5, -.035) [scale=1, color=black] {$\frac{\sf 3}{\sf 5}$};
\node at (4/5, -.035) [scale=1, color=black] {$\frac{\sf 4}{\sf 5}$};
\node at (1/6, -.035) [scale=1, color=black] {$\frac{\sf 1}{\sf 6}$};

%\node at (0/3, 2/3) [scale=.9, color=black] {\sf b};
%\node at (1/2, 0/2) [scale=.9, color=black] {\sf c};

\end{tikzpicture}
\caption{Ford circles}
\label{fig:Ford}
\end{figure}

\begin{proposition}
\label{thm:Ford2}
\sf
The point of tangency of two tangent Ford circles \eqref{eq:matrix} is 
%$$
\begin{equation}
\label{eq:point}
%\det\begin{bmatrix}p&q\\m&n\end{bmatrix} = \pm1
%\qquad \Rightarrow\qquad
K[p,m] \cap K[q,n] \ = \  \left(\frac{pm+qn}{m^2+n^2},\, \frac{1}{m^2+n^2} \right)
\end{equation}
%$$
\end{proposition}

\noindent
{\bf Proof:}  
Simple similarity of triangles leads to $(\frac{q}{n} - x) : (x-\frac{p}{m}) = R:r$
with the radii being $R=\frac{1}{2n^2}$ and $r=\frac{1}{2m^2}$,
which easily solves for $x$.
\qed

~

\noindent
{\bf Remark:}
Remark:  Note that formula \eqref{eq:point} differs from the one given in \cite{Apostol} 
and reported in \cite{Weisstein};
it is symmetric and does not depend on the order of circles.

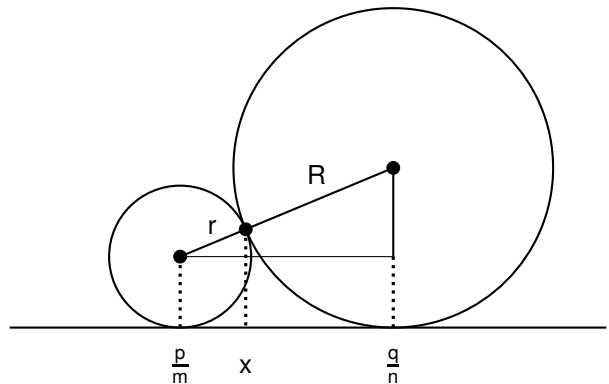
\begin{figure}[H]
\centering
\begin{tikzpicture}[scale=17]
%\clip (-.1, -.12) rectangle (1.1, .56);
%%-----------------------  
\draw [thick] (1/5, 0)--(1/2 + 1/6,0);  %line

\foreach \a /  \b  in {
1/2, 1/3
}
\draw [thick]  (\a/\b, 1/2*1/\b*1/\b) circle (1/2*1/\b*1/\b);

\draw [fill=black] (5/13,1/13) circle (1/200); %
\draw [fill=black] (1/3,1/18) circle (1/200); %
\draw [fill=black] (1/2,1/8) circle (1/200); %diagonal

\draw [-] (1/3, 1/18) -- (1/2, 1/18);
\draw [-,thick] (1/3, 1/18) -- (1/2, 1/8); 
\draw [-,thick] (1/2, 1/18) -- (1/2, 1/8); 

\draw [very thick, dotted] (1/3, 1/18) -- (1/3, 0); 
\draw [very thick, dotted] (1/2, 1/18) -- (1/2, 0); 
\draw [very thick, dotted] (5/13, 1/13) -- (5/13, 0); 

\node at (1/3+.025, .079) [scale=1, color=black] {\sf r};
\node at (1/2-.06, .12) [scale=1, color=black] {\sf R};

\node at (1/2, -.03) [scale=1, color=black] {$\frac{\sf q}{\sf n}$};
\node at (1/3, -.03) [scale=1, color=black] {$\frac{\sf p}{\sf m}$};
\node at (5/13, -.03) [scale=1, color=black] {${\sf x}$};

\end{tikzpicture}
\caption{Derivation of formula \eqref{eq:point}}
\label{fig:proof}
\end{figure}

\newpage

%-------------------------------------------------------------
\section{Fibonacci numbers}

By the extended Fibonacci sequence we understand a bilateral sequence 
$F_n$ defined by:
$$
F_0=0, \quad F_1=1,\quad F_{n+1}= F_n+F_{n-1}, \quad n\in \mathbb Z
$$
It looks like this
$$
\begin{array}{ccccccccccccccccccccccccccccccccccccccccccccc} 
                    ...   
                &   5    
                &  -3  
                &   2    
                &  -1   
                &   1    
                &  0   
                &  1  
                &  1  
                &  2  
                &  3                  
                &  5                 
                &  8             
                &  13 
                &  21                  
                &  34
                & ...
                \\[7pt]
                   
                &     
                &    
                & \uparrow   
                &  \uparrow   
                &  \uparrow   
                & \uparrow   
                &  \uparrow   
                &  \uparrow   
                &  \uparrow   
                &  \uparrow   
                &  \uparrow   
                &              
                &   
                &                   
                &  
                &                    
                       \\
                &     
                &    
                &  \mbox{\footnotesize $ F_{-3} $ }   
                &  \mbox{\footnotesize $ F_{-2} $ }   
                &  \mbox{\footnotesize $ F_{-1} $ }   
                &  \mbox{\footnotesize $ F_{0} $ }   
                &  \mbox{\footnotesize $ F_{1} $ }   
                &  \mbox{\footnotesize $ F_{2} $ }   
                &  \mbox{\footnotesize $ F_{3} $ }   
                &  \mbox{\footnotesize $ F_{4} $ }   
                &  \mbox{\footnotesize $ F_{5} $ }   
                &              
                &   
                &                   
                &  
                & 
\end{array}
$$
The extended Fibonacci sequence satisfies the basic identities valid for the standard Fibonacci sequence.
Among them are the following:
%$$
\begin{equation}
\label{eq:identities}
\begin{array}{cllll}
(i) & F_{n-1}F_{n}+F_{n}F_{n+1} \ = \  F_{2n}\\[9pt]
(ii) & F_{n}^2 +F_{n+1}^2 \ = \ F_{2n+1}\\[9pt]
(iii) & F_{n+1}F_{n-1} - F_n^2 \ = \  (-1)^n
\end{array}
\end{equation}
%$$

%
%\, \frac{1}{F_{n}^2 +F_{n+1}^2} \right)
%     \ = \ 
%     \left(\frac{}{F_{2n+1}},\, \frac{1}} \right)
%
%

%\newpage

%------------------------------------------------------------------
\section{Ford and Fibonacci}

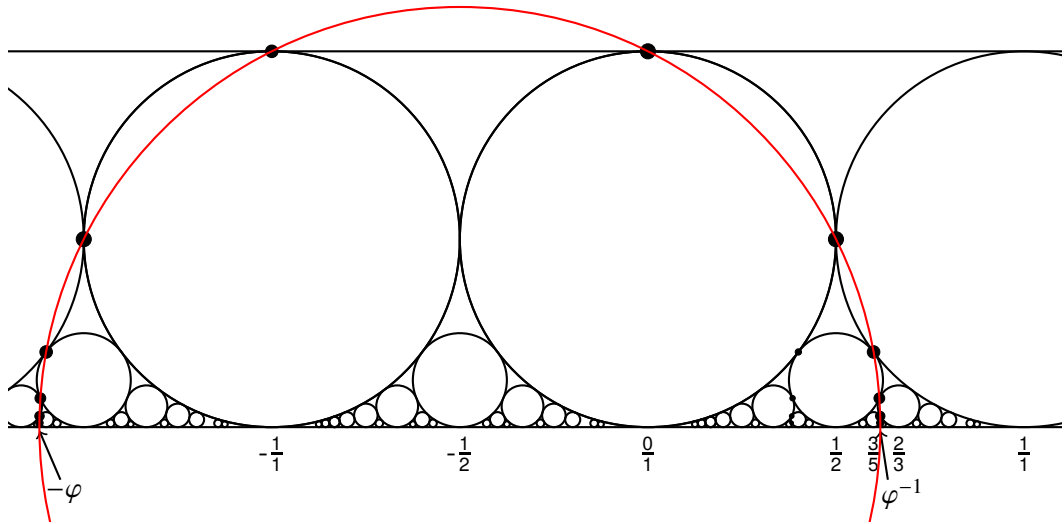
\begin{tikzpicture}[scale=5]
\clip (-1.7, -.25) rectangle (1.1, 1.2);
%%-----------------------  
\draw [thick] (-3,0)--(2,0);  %line
\draw [thick] (-3,1)--(2,1);  %line

\foreach \a /  \b  in {
0/1, 1/1,
1/2,
1/3,2/3,
1/4,3/4,
1/5,2/5,3/5,4/5,
1/6,1/6,
1/7,2/7,3/7,4/7,5/7,6/7,
1/8,3/8,5/8,7/8
}
\draw [thick]  (\a/\b, 1/2*1/\b*1/\b) circle (1/2*1/\b*1/\b)
%   (1+\a/\b, 1/2*1/\b*1/\b) circle (1/2*1/\b*1/\b)
   (-1+\a/\b, 1/2*1/\b*1/\b) circle (1/2*1/\b*1/\b)
   (-2+\a/\b, 1/2*1/\b*1/\b) circle (1/2*1/\b*1/\b);

%\draw  [->,  red, bend right=-20] (1/2,1/2) to  (3/5, 1/5);

\draw [fill=black] (1/2,1/2) circle (2/100); %diagonal
\draw [fill=black] (3/5,1/5) circle (2/120); %diagonal
\draw [fill=black] (8/13,1/13) circle (2/140); %diagonal
\draw [fill=black] (21/34,1/34) circle (2/160); %diagonal
\draw [fill=black] (55/89,1/89) circle (2/200); %diagonal

\draw [fill=black] (0,1) circle (2/100); %diagonal
\draw [fill=black] (-1,1) circle (2/120); %diagonal

\draw [fill=black] (-1-1/2,1/2) circle (2/100); %diagonal
\draw [fill=black] (-1-3/5,1/5) circle (2/120); %diagonal
\draw [fill=black] (-1-8/13,1/13) circle (2/140); %diagonal
\draw [fill=black] (-1-21/34,1/34) circle (2/160); %diagonal
\draw [fill=black] (-1-55/89,1/89) circle (2/200); %diagonal

\draw [fill=black] (1/2,1/2) circle (1/100); %diagonal
\draw [fill=black] (2/5,1/5) circle (1/120); %diagonal
\draw [fill=black] (5/13,1/13) circle (1/140); %diagonal
\draw [fill=black] (13/34,1/34) circle (1/160); %diagonal
\draw [fill=black] (34/89,1/89) circle (1/200); %diagonal

\node at (0/1, -.07) [scale=1, color=black] {$\frac{\sf 0}{\sf 1}$};
\node at (1/1, -.07) [scale=1, color=black] {$\frac{\sf 1}{\sf 1}$};
\node at (1/2, -.07) [scale=1, color=black] {$\frac{\sf 1}{\sf 2}$};
%\node at (1/3, -.07) [scale=1, color=black] {$\frac{\sf 1}{\sf 3}$};
\node at (2/3, -.07) [scale=1, color=black] {$\frac{\sf 2}{\sf 3}$};
%\node at (1/4, -.07) [scale=1, color=black] {$\frac{\sf 1}{\sf 4}$};
%\node at (3/4, -.07) [scale=1, color=black] {$\frac{\sf 3}{\sf 4}$};
%\node at (1/5, -.07) [scale=1, color=black] {$\frac{\sf 1}{\sf 5}$};
%\node at (2/5, -.07) [scale=1, color=black] {$\frac{\sf 2}{\sf 5}$};
\node at (3/5, -.07) [scale=1, color=black] {$\frac{\sf 3}{\sf 5}$};
%\node at (4/5, -.07) [scale=1, color=black] {$\frac{\sf 4}{\sf 5}$};
%\node at (1/6, -.07) [scale=1, color=black] {$\frac{\sf 1}{\sf 6}$};

\node at (-1/1, -.07) [scale=1, color=black] {-$\frac{\sf 1}{\sf 1}$};
\node at (-1/2, -.07) [scale=1, color=black] {-$\frac{\sf 1}{\sf 2}$};

%\draw [->,thick] (.382, -.14) -- (.382, -.01);
\draw [->,thick] (.618+.02, -.14) -- (.618, -.01); 
\draw [->,thick] (-1.618+.056, -.14) -- (-1.618, -.01);

%\node at (1/2.618, -.177) [scale=1, color=black] {$\ \varphi^{-2}$};
\node at (1/1.618+.049, -.177) [scale=1, color=black] {$\ \varphi^{-1}$};
\node at (-1.618+.056, -.177) [scale=1, color=black] {$\ -\varphi$};

%\node at (0/3, 2/3) [scale=.9, color=black] {\sf b};
%\node at (1/2, 0/2) [scale=.9, color=black] {\sf c};

\draw [thick, red] (-1/2,0) circle (1/2+.618); %circle

%\draw [thick, blue] (3/2,0) circle (1/2+.618); %circle

\end{tikzpicture}

\begin{proposition}
\label{thm:Ford3}
\sf
Let $F_n$ be the extended Fibonacci sequence:
%$$
%F_0=0, \quad F_1=1,\quad F_{n+1}= F_n+F_{n-1}, \quad n\in \mathbb Z
%$$
Then the sequence of points with coordinates 
%$$
\begin{equation}
\label{eq:A}
\left( \frac{F_{2n}}{F_{2n+1}},\,\frac{1}{F_{2n+1}}\right)
% \qquad   (red)
%...,\, 
%\left( \frac{-8}{5},\,\frac{1}{5}\right),\,
%\left( \frac{-3}{2},\,\frac{1}{2}\right),\,
%\left( \frac{-1}{1},\,\frac{1}{1}\right),\,
%\left( \frac{0}{1},\,\frac{1}{1}\right),\,
%\left( \frac{1}{2},\,\frac{1}{2}\right),\,
%\left( \frac{3}{5},\,\frac{1}{5}\right),\,
%\left( \frac{8}{13},\,\frac{1}{13}\right),\,
%\left( \frac{21}{34},\,\frac{1}{34}\right)...
\end{equation}
%$$
lie on a circle centered at  $(-\frac{1}{2}, 0)$ of radius $\sqrt{5}/2$.
Moreover, the above points are among points of tangency of Ford circles.
The circle cuts the x-axis at the golden cut $x=\varphi^{-1}$ and 
negative golden ratio, $x=-\varphi$.
\end{proposition}

\noindent
{\bf Proof:}
First, let us prove the the points are concyclic.
Denote $2n=m$. We need to verify the quadratic circle equation 
$$
\left(\frac{F_m}{F_{m+1}} + \frac{1}{2}\right)^2 + \left(\frac{1}{F_{m+1}}\right)^2 = \frac{5}{4}\,.
$$
After squaring, collecting like terms, and factoring out the common 4, we get
$$
F_{m+1}^2 - F_mF_{m+1} =F_m^2 +1 \,.
$$
Next, factor the common $F_{m+1}$:
$$
F_{m+1}\, (F_{m+1} - F_m) = F_m^2+1
$$
and use the Fibonacci's defining property to get:
$$
F_{m+1}F_{m-1} = F_m^2 +1 \,,
$$
which is one of the standard identities true for even $m$,
see Eq. \eqref{eq:identities} ({\it iii}).
\\

For the second part,
consider the Ford circles 
above the Fibonacci fractions:
$$
\frac{F_n}{F_{n+1}}:\qquad 
...,\,
\frac{2}{-1},\,
\frac{-1}{1},\,
\frac{1}{0},\,
\frac{0}{1},\,
\frac{1}{1},\,
\frac{1}{2},\,
\frac{2}{3},\,
\frac{3}{5},\,
\frac{5}{8},\,
...
$$
The consecutive Ford circles over this sequence are clearly mutually tangent, since 
$$
\det\begin{bmatrix}\;F_{n-1}&F_n\\F_n&F_{n+1}\;\end{bmatrix} = (-1)^{n}
$$
thus \eqref{eq:matrix} is satisfied.
What remains is to check what are their points of tangency.
For this, we use Proposition \ref{thm:Ford2}.  Under substitution we get:
%$$
%K[F_n,F_{n+1}] \cap K[F_{n+1}{F_{n+2}}] \ = \
%     \left(\frac{F_nF_{n+1}+F_{n+1}F_{n+2}}{F_{n+1}^2 +F_{n+2}^2},\, \frac{1}{F_{n+1}^2 +F_{n+2}^2} \right)
%     \ = \ 
%     \left(\frac{F_{2n+2}}{F_{2n+3}},\, \frac{1}{F_{2n+3}} \right)
%$$
$$
K[F_{n-1},F_{n}] \cap K[F_{n},{F_{n+1}}] \ = \
     \left(\frac{F_{n-1}F_{n}+F_{n}F_{n+1}}{F_{n}^2 +F_{n+1}^2},\, \frac{1}{F_{n}^2 +F_{n+1}^2} \right)
     \ = \ 
     \left(\frac{F_{2n}}{F_{2n+1}},\, \frac{1}{F_{2n+1}} \right)
$$
where we used \eqref{eq:identities} ({\it i}) and ({\it ii}\/).
This concludes the proof.
\qed

\newpage

%-----------------------------------------------------------------------------------
\section{Four concyclic sequences}

There are four instances of similar correspondence between Fibonacci sequence and Ford circles.
Here they are denoted as (A), (B), (C), and (D), together with the defining recipe:

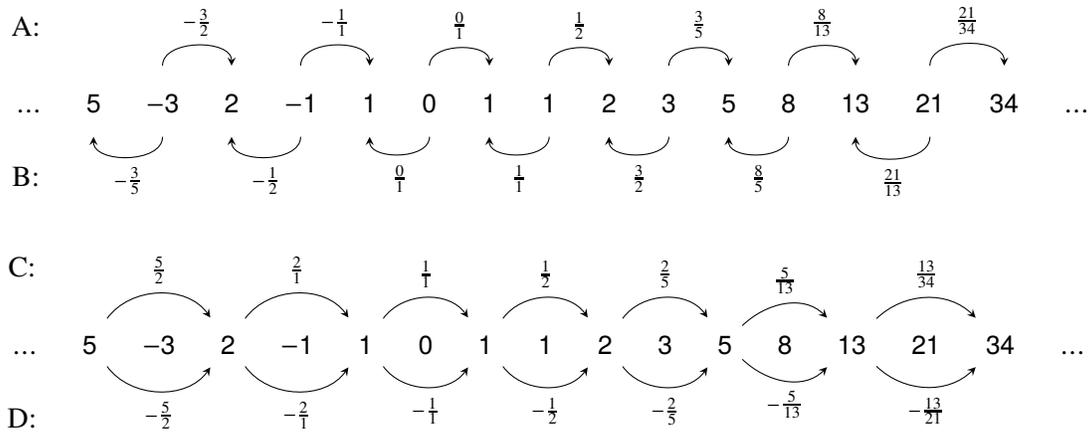
\hspace{-.48in}
\begin{tikzpicture}[baseline=-0.7ex]
    \matrix (m) [ matrix of math nodes,
                         row sep=1.5em,
                         column sep=.9em,
                         text height=2.8ex, text depth=1.1ex] 
      {
                    ...   
                &  \sf  5    
                &  \sf -3  
                &  \sf  2    
                &  \sf -1   
                &  \sf  1    
                &  \sf 0   
                &  \sf 1  
                &  \sf 1  
                &  \sf 2  
                &  \sf 3                  
                &  \sf 5                 
                &  \sf 8             
                &  \sf 13 
                &  \sf 21                  
                &  \sf 34
                & ...
                \\
     };

    \path[-stealth]
%horizontal
        (m-1-3) edge [out=90, in=90]  node[above ] { \smalll $-\frac{3}{2}$}  (m-1-4)
        (m-1-5) edge [out=90, in=90]  node[above ] { \smalll $-\frac{1}{1}$}  (m-1-6)
        (m-1-7) edge [out=90, in=90]  node[above ] { \smalll $\frac{0}{1}$}  (m-1-8)
        (m-1-9) edge [out=90, in=90]  node[above ] { \smalll $\frac{1}{2}$}  (m-1-10)        
        (m-1-11) edge [out=90, in=90]  node[above ] { \smalll $\frac{3}{5}$}  (m-1-12)        
        (m-1-13) edge [out=90, in=90]  node[above ] { \smalll $\frac{8}{13}$}  (m-1-14)
         (m-1-15) edge [out=90, in=90]  node[above ] { \smalll $\frac{21}{34}$}  (m-1-16)
        
        (m-1-3) edge [out=-90, in=-90]  node[below ] { \smalll $-\frac{3}{5}$}  (m-1-2)
        (m-1-5) edge [out=-90, in=-90]  node[below ] { \smalll $-\frac{1}{2}$}  (m-1-4)
        (m-1-7) edge [out=-90, in=-90]  node[below ] { \smalll $\frac{0}{1}$}  (m-1-6)
        (m-1-9) edge [out=-90, in=-90]  node[below ] { \smalll $\frac{1}{1}$}  (m-1-8)
        (m-1-11) edge [out=-90, in=-90]  node[below ] { \smalll $\frac{3}{2}$}  (m-1-10)
        (m-1-13) edge [out=-90, in=-90]  node[below ] { \smalll $\frac{8}{5}$}  (m-1-12)
        (m-1-15) edge [out=-90, in=-90]  node[below ] { \smalll $\frac{21}{13}$}  (m-1-14)
;
\node at (-7,1) {A:};
\node at (-7,-1) {B:};
\end{tikzpicture}   
%\end{equation}

~\\

\hspace{-.5in}
\begin{tikzpicture}[baseline=-0.7ex]
    \matrix (m) [ matrix of math nodes,
                         row sep=1.5em,
                         column sep=.9em,
                         text height=2.8ex, text depth=1.1ex] 
      {
                    ...   
                &   \sf 5    
                &  \sf -3  
                &  \sf  2    
                &  \sf -1   
                &  \sf  1    
                &  \sf 0   
                &  \sf 1  
                &  \sf 1  
                &  \sf 2  
                &  \sf 3                  
                &  \sf 5                 
                &  \sf 8             
                &  \sf 13 
                &  \sf 21                  
                &  \sf 34
                & ...
                \\
     };

    \path[-stealth]
%horizontal
        (m-1-2) edge [out=54, in=125]  node[above ] { \smalll $\frac{5}{2}$}  (m-1-4)
        (m-1-4) edge [out=54, in=125]  node[above ] { \smalll $\frac{2}{1}$}  (m-1-6)
        (m-1-6) edge [out=54, in=125]  node[above ] { \smalll $\frac{1}{1}$}  (m-1-8)
        (m-1-8) edge [out=54, in=125]  node[above ] { \smalll $\frac{1}{2}$}  (m-1-10)        
        (m-1-10) edge [out=54, in=125]  node[above ] { \smalll $\frac{2}{5}$}  (m-1-12)        
        (m-1-12) edge [out=44, in=135]  node[above ] { \smalll $\frac{5}{13}$}  (m-1-14)
         (m-1-14) edge [out=44, in=135]  node[above ] { \smalll $\frac{13}{34}$}  (m-1-16)
        
        (m-1-2) edge [out=-54, in=-125]  node[below ] { \smalll $-\frac{5}{2}$}  (m-1-4)
        (m-1-4) edge [out=-54, in=-125]  node[below ] { \smalll $-\frac{2}{1}$}  (m-1-6)
        (m-1-6) edge [out=-54, in=-125]  node[below ] { \smalll $-\frac{1}{1}$}  (m-1-8)
        (m-1-8) edge [out=-54, in=-125]  node[below ] { \smalll $-\frac{1}{2}$}  (m-1-10)
        (m-1-10) edge [out=-54, in=-125]  node[below ] { \smalll $-\frac{2}{5}$}  (m-1-12)
        (m-1-12) edge [out=-44, in=-135]  node[below ] { \smalll $-\frac{5}{13}$}  (m-1-14)
        (m-1-14) edge [out=-44, in=-135]  node[below ] { \smalll $-\frac{13}{21}$}  (m-1-16)
;
\node at (-7,1) {C:};
\node at (-7,-1) {D:};
\end{tikzpicture}   
%\end{equation}

Sequence (A) is the one discussed in the previous sections.
Each sequence corespondents to a circle of radius $\sqrt{5}/2$ and intersects ford circles at points 
the coordinates of which are ratios of Fibonacci numbers.
The centers lie on the x-axis at the following points: 
$$ 
A:  \ x=-1/2,\ \   
B:  \  x= 1/2, \ \   
C:  \  x= 3/2, \  \ %(blue, symmetric, step 2)
D: \  x=-3/2  %(yellow, step 2, symmetric) 
 $$

~

%$$
%\left(\frac{1}{2},\, \frac{1}{2} \right), \quad
%\left(\frac{3}{5},\, \frac{1}{5} \right), \quad
%\left(\frac{8}{13},\, \frac{1}{13} \right), \quad
%...
%\quad
%\rightarrow \  \frac{1}{\varphi} = \frac{\sqrt{5}-1}{2}
%$$
%
%left:
%$$
%\left(\frac{1}{2},\, \frac{1}{2} \right), \quad
%\left(\frac{2}{5},\, \frac{1}{5} \right), \quad
%\left(\frac{5}{13},\, \frac{1}{13} \right), \quad
%...
%\quad
%\rightarrow \  \frac{1}{\varphi^2} = \frac{3-\sqrt{5}}{2}
%$$
%

\begin{figure}\hspace{-.35in}
\begin{tikzpicture}[scale=2.7]
\clip (-2.7, -.37) rectangle (2.7, 1.28);
%%-----------------------  
\draw [thick] (-3,0)--(3,0);  %line
\draw [thick] (-3,1)--(3,1);  %line

\foreach \a /  \b  in {
0/1, 1/1,
1/2,
1/3,2/3,
1/4,3/4,
1/5,2/5,3/5,4/5,
1/6,1/6,
1/7,2/7,3/7,4/7,5/7,6/7,
1/8,3/8,5/8,7/8
}
\draw [thick]  (\a/\b, 1/2*1/\b*1/\b) circle (1/2*1/\b*1/\b)
   (2+\a/\b, 1/2*1/\b*1/\b) circle (1/2*1/\b*1/\b)
   (1+\a/\b, 1/2*1/\b*1/\b) circle (1/2*1/\b*1/\b)
   (-1+\a/\b, 1/2*1/\b*1/\b) circle (1/2*1/\b*1/\b)
   (-2+\a/\b, 1/2*1/\b*1/\b) circle (1/2*1/\b*1/\b)
   (-3+\a/\b, 1/2*1/\b*1/\b) circle (1/2*1/\b*1/\b)
;

%\draw  [->,  red, bend right=-20] (1/2,1/2) to  (3/5, 1/5);

\draw [fill=black] (1/2,1/2) circle (2/100); %diagonal
\draw [fill=black] (3/5,1/5) circle (2/120); %diagonal
\draw [fill=black] (8/13,1/13) circle (2/140); %diagonal
\draw [fill=black] (21/34,1/34) circle (2/160); %diagonal
\draw [fill=black] (55/89,1/89) circle (2/200); %diagonal

\draw [fill=black] (0,1) circle (2/100); %diagonal
\draw [fill=black] (-1,1) circle (2/120); %diagonal

\draw [fill=black] (-1-1/2,1/2) circle (2/100); %diagonal
\draw [fill=black] (-1-3/5,1/5) circle (2/120); %diagonal
\draw [fill=black] (-1-8/13,1/13) circle (2/140); %diagonal
\draw [fill=black] (-1-21/34,1/34) circle (2/160); %diagonal
\draw [fill=black] (-1-55/89,1/89) circle (2/200); %diagonal

\draw [fill=black] (1/2,1/2) circle (1/100); %diagonal
\draw [fill=black] (2/5,1/5) circle (1/120); %diagonal
\draw [fill=black] (5/13,1/13) circle (1/140); %diagonal
\draw [fill=black] (13/34,1/34) circle (1/160); %diagonal
\draw [fill=black] (34/89,1/89) circle (1/200); %diagonal

\node at (-3, -.11) [scale=1, color=black] {$\sf -3$};
\node at (-2, -.11) [scale=1, color=black] {$\sf -2$};
\node at (-1, -.11) [scale=1, color=black] {$\sf -1$};
\node at (0/1, -.11) [scale=1, color=black] {$0$};
\node at (1, -.11) [scale=1, color=black] {$\sf 1$};
\node at (2, -.11) [scale=1, color=black] {$\sf 2$};
\node at (3, -.11) [scale=1, color=black] {$\sf 3$};
%\node at (1/1, -.11) [scale=1, color=black] {$\sf 1$};

%
%\node at (-1/2, -.11) [scale=1, color=black] {$\frac{\sf 1}{\sf 2}$};
%\node at (1/2, -.11) [scale=1, color=black] {$\frac{\sf 1}{\sf 2}$};
%\node at (3/2, -.11) [scale=1, color=black] {$\frac{\sf 1}{\sf 2}$};
%\node at (5/2, -.11) [scale=1, color=black] {$\frac{\sf 1}{\sf 2}$};
%
%\node at (1/3, -.07) [scale=1, color=black] {$\frac{\sf 1}{\sf 3}$};
%\node at (2/3, -.07) [scale=1, color=black] {$\frac{\sf 2}{\sf 3}$};
%\node at (1/4, -.07) [scale=1, color=black] {$\frac{\sf 1}{\sf 4}$};
%\node at (3/4, -.07) [scale=1, color=black] {$\frac{\sf 3}{\sf 4}$};
%\node at (1/5, -.07) [scale=1, color=black] {$\frac{\sf 1}{\sf 5}$};
%\node at (2/5, -.07) [scale=1, color=black] {$\frac{\sf 2}{\sf 5}$};
%\node at (3/5, -.07) [scale=1, color=black] {$\frac{\sf 3}{\sf 5}$};
%\node at (4/5, -.07) [scale=1, color=black] {$\frac{\sf 4}{\sf 5}$};
%\node at (1/6, -.07) [scale=1, color=black] {$\frac{\sf 1}{\sf 6}$};

\draw [->,thick] (-2.62, -.17) -- (-2.62, -.05);
%\draw [->,thick] (-2.79, -.21) -- (-2.65, -.05); 
\node at (-2.62, -.3) [scale=1, color=black] {$\ -\varphi^2$};
%\node at (-2.78, -.28) [scale=1, color=black] {-$\ \varphi$};

%\draw [->,thick] (-1.21, -.21) -- (-1.35, -.05);
\draw [->,thick] (-1.62, -.17) -- (-1.62, -.05); 
%\node at (-1.22, -.28) [scale=1, color=black] {-$\ \varphi^{-2}$};
\node at (-1.62, -.3) [scale=1, color=black] {$\ -\varphi\phantom{^{2}}$};

\draw [->,thick] (-.38, -.17) -- (-.38, -.05);
\draw [->,thick] (-.62, -.17) -- (-.62, -.05); 
\node at (-.38, -.3) [scale=1, color=black] {$-\tau^2$};
\node at (-.62, -.3) [scale=1, color=black] {$-\tau\phantom{^{2}}\ $};

\draw [->,thick] (.38, -.17) -- (.38, -.05);
\draw [->,thick] (.62, -.17) -- (.62, -.05); 
\node at (.38, -.3) [scale=1, color=black] {$\ \tau^{2}$};
\node at (.62, -.3) [scale=1, color=black] {$\ \tau\phantom{^{2}}$};

%\draw [->,thick] (1.21, -.21) -- (1.35, -.05);
\draw [->,thick] (1.62, -.17) -- (1.62, -.05); 
%\node at (1.22, -.28) [scale=1, color=black] {$\ \varphi^{-2}$};
\node at (1.62, -.3) [scale=1, color=black] {$\ \varphi\phantom{^{2}}$};

\draw [->,thick] (2.62, -.17) -- (2.62, -.05);
%\draw [->,thick] (2.79, -.21) -- (2.65, -.05); 
\node at (2.62, -.3) [scale=1, color=black] {$\ \varphi^{2}$};
%\node at (2.78, -.28) [scale=1, color=black] {$\ \varphi^2$};

\draw [-stealth, very thick] (0,0) -- (0,1.28);

%\node at (0/3, 2/3) [scale=.9, color=black] {\sf b};
%\node at (1/2, 0/2) [scale=.9, color=black] {\sf c};

\begin{scope}
\clip (-3,0) rectangle (3,2); 
\draw [thick, red] (-1/2,0) circle (1/2+.618); %circle
\draw [thick, blue] (3/2,0) circle (1/2+.618); %circle

\draw [thick, red] (1/2,0) circle (1/2+.618); %circle
\draw [thick, blue] (-3/2,0) circle (1/2+.618); %circle
\end{scope}

\node at (-3/2, 1.2) [scale=1, color=black] {\sf D};
\node at (-1/2, 1.2) [scale=1, color=black] {\sf A};
\node at (1/2, 1.2) [scale=1, color=black] {\sf B};
\node at (3/2, 1.2) [scale=1, color=black] {\sf C};

\end{tikzpicture}
\caption{Four golden circles upon Ford circles}
\label{fig:four}
\end{figure}
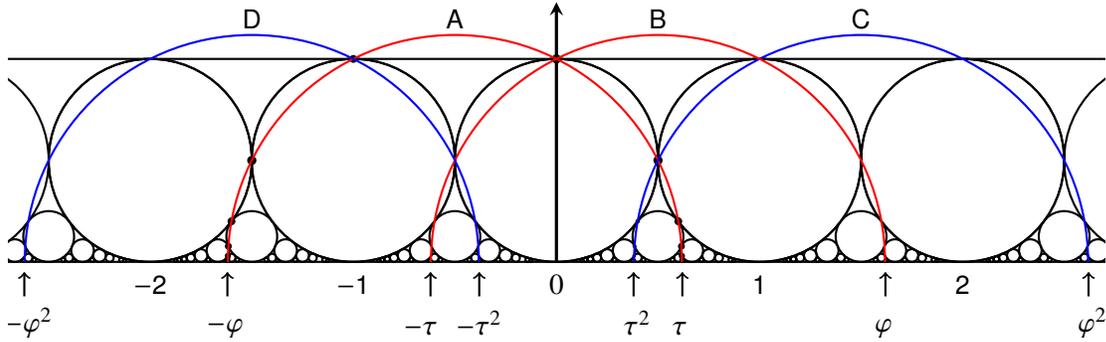

The circles cut the $x$-axis at different values of the golden ratio,
as show in Figure \ref{fig:four}.
We use the following notation:
$$
\varphi = \frac{\sqrt{5}+1}{2} \approx 1.618\,,
\qquad
\tau = \frac{\sqrt{5}-1}{2} \approx 0.618\,,
\qquad 
\varphi \tau = 1
$$

~\\
{\bf Proof:}
One may start with sequence (A) and add 1 to the $x$variable in  every term:
$$
%\begin{equation}
%\label{eq:B}
\left( \frac{F_{2n}}{F_{2n+1}}+1,\,\frac{1}{F_{2n+1}}\right)
\ = \ 
\left( \frac{F_{2n+2}}{F_{2n+1}},\,\frac{1}{F_{2n+1}}\right)
% \qquad   (red)
%...,\, 
%\left( \frac{-8}{5},\,\frac{1}{5}\right),\,
%\left( \frac{-3}{2},\,\frac{1}{2}\right),\,
%\left( \frac{-1}{1},\,\frac{1}{1}\right),\,
%\left( \frac{0}{1},\,\frac{1}{1}\right),\,
%\left( \frac{1}{2},\,\frac{1}{2}\right),\,
%\left( \frac{3}{5},\,\frac{1}{5}\right),\,
%\left( \frac{8}{13},\,\frac{1}{13}\right),\,
%\left( \frac{21}{34},\,\frac{1}{34}\right)...
%\end{equation}
$$
which leads to sequence (B).
Shifting again by 1 produces
$$
%\begin{equation}
%\label{eq:C}
\left( \frac{F_{2n+2}}{F_{2n+1}}+1,\,\frac{1}{F_{2n+1}}\right)
\ = \ 
\left( \frac{F_{2n+3}}{F_{2n+1}},\,\frac{1}{F_{2n+1}}\right)
% \qquad   (red)
%...,\, 
%\left( \frac{-8}{5},\,\frac{1}{5}\right),\,
%\left( \frac{-3}{2},\,\frac{1}{2}\right),\,
%\left( \frac{-1}{1},\,\frac{1}{1}\right),\,
%\left( \frac{0}{1},\,\frac{1}{1}\right),\,
%\left( \frac{1}{2},\,\frac{1}{2}\right),\,
%\left( \frac{3}{5},\,\frac{1}{5}\right),\,
%\left( \frac{8}{13},\,\frac{1}{13}\right),\,
%\left( \frac{21}{34},\,\frac{1}{34}\right)...
%\end{equation}
$$
which defines the sequence and circle (C). 
Shifting (A) by a unit to the left produces
$$
%\begin{equation}
%\label{eq:B}
\left( \frac{F_{2n}}{F_{2n+1}}-1,\,\frac{1}{F_{2n+1}}\right)
\ = \ 
\left( \frac{-F_{2n-1}}{F_{2n+1}},\,\frac{1}{F_{2n+1}}\right)
% \qquad   (red)
%...,\, 
%\left( \frac{-8}{5},\,\frac{1}{5}\right),\,
%\left( \frac{-3}{2},\,\frac{1}{2}\right),\,
%\left( \frac{-1}{1},\,\frac{1}{1}\right),\,
%\left( \frac{0}{1},\,\frac{1}{1}\right),\,
%\left( \frac{1}{2},\,\frac{1}{2}\right),\,
%\left( \frac{3}{5},\,\frac{1}{5}\right),\,
%\left( \frac{8}{13},\,\frac{1}{13}\right),\,
%\left( \frac{21}{34},\,\frac{1}{34}\right)...
%\end{equation}
$$
which produces thesequence  and circle (D).
\qed

~

Sequences (A) (B) and (C) are genuine ratios of the Fibonacci numbers, while 
sequence (D) has negative signs of the $x$-coefficients imposed artificially.
Note that A and B are negatives of each other, and so are C and D.
The interval (0,1) contains only fragments of (A) and (B).

~

\noindent
{\bf Remark:}  
Note that the open circle $y=1$ contributes to the points in the sequences (A), (B), (C), and (D).

%\noindent
%E.g., the sequence B consists of points 
%with coordinates 
%$$
%\left( \frac{F_{2n-1}}{F_{2n+1}},\,\frac{1}{F_{2n+1}}\right)
%\qquad (blue)
%$$
%lie on a circle centered at  $(\frac{3}{2}, 0)$ of radius $\sqrt{5}/2$.
%The above points are among points of tangency of Ford circles as well.
%

~

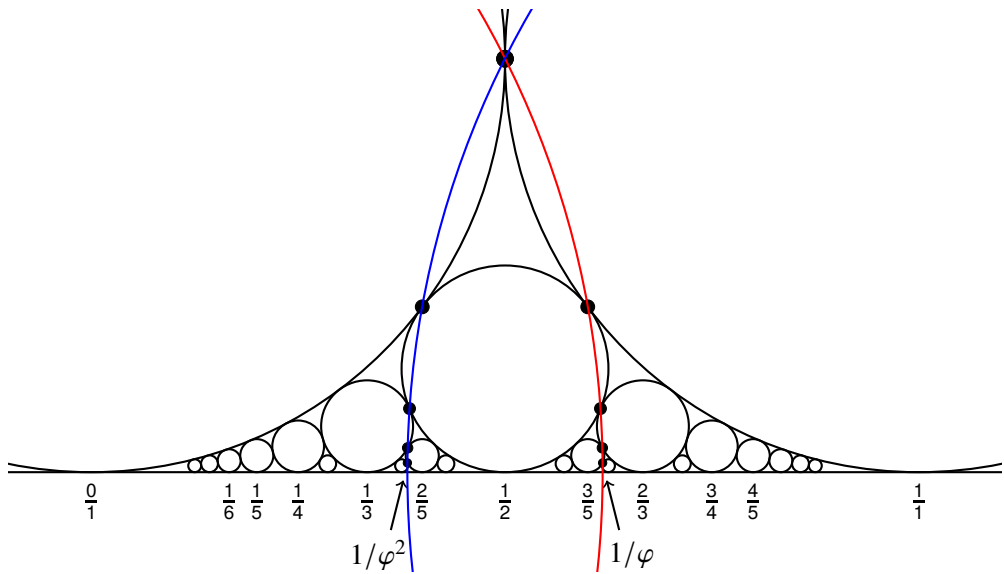
\begin{figure}[H]
\begin{tikzpicture}[scale=11]
\clip (-.1, -.12) rectangle (1.1, .56);
%%-----------------------  
\draw [thick] (-.2,0)--(1.2,0);  %line

\foreach \a /  \b  in {
0/1, 1/1, 1/2, 1/3,2/3, 1/4,3/4, 1/5,2/5,3/5,4/5,
1/6,5/6, 1/7,2/7,3/7,4/7,5/7,6/7, 1/8,3/8,5/8,7/8
}
\draw [thick]  (\a/\b, 1/2*1/\b*1/\b) circle (1/2*1/\b*1/\b)
%   (1+\a/\b, 1/2*1/\b*1/\b) circle (1/2*1/\b*1/\b)
%   (-1+\a/\b, 1/2*1/\b*1/\b) circle (1/2*1/\b*1/\b)
%  (-2+\a/\b, 1/2*1/\b*1/\b) circle (1/2*1/\b*1/\b)
;

\draw [fill=black] (1/2,1/2) circle (1/100); %diagonal
\draw [fill=black] (3/5,1/5) circle (1/120); %diagonal
\draw [fill=black] (8/13,1/13) circle (1/140); %diagonal
\draw [fill=black] (21/34,1/34) circle (1/160); %diagonal
\draw [fill=black] (55/89,1/89) circle (1/200); %diagonal

%\draw  [->,  red, bend right=-20] (1/2,1/2) to  (3/5, 1/5);

\draw [fill=black] (1/2,1/2) circle (1/100); %diagonal
\draw [fill=black] (2/5,1/5) circle (1/120); %diagonal
\draw [fill=black] (5/13,1/13) circle (1/140); %diagonal
\draw [fill=black] (13/34,1/34) circle (1/160); %diagonal
\draw [fill=black] (34/89,1/89) circle (1/200); %diagonal

\draw [->,thick] (.382-.02, -.07) -- (.382-.005, -.01);
\draw [->,thick] (.618+.02, -.07) -- (.618+.005, -.01); 

\node at (1/2.618-.035, -.1) [scale=1, color=black] {$1/\varphi^2$};
\node at (1/1.618+.035, -.1) [scale=1, color=black] {$1/\varphi$};

\node at (0/1, -.035) [scale=1, color=black] {$\frac{\sf 0}{\sf 1}$};
\node at (1/1, -.035) [scale=1, color=black] {$\frac{\sf 1}{\sf 1}$};
\node at (1/2, -.035) [scale=1, color=black] {$\frac{\sf 1}{\sf 2}$};
\node at (1/3, -.035) [scale=1, color=black] {$\frac{\sf 1}{\sf 3}$};
\node at (2/3, -.035) [scale=1, color=black] {$\frac{\sf 2}{\sf 3}$};
\node at (1/4, -.035) [scale=1, color=black] {$\frac{\sf 1}{\sf 4}$};
\node at (3/4, -.035) [scale=1, color=black] {$\frac{\sf 3}{\sf 4}$};
\node at (1/5, -.035) [scale=1, color=black] {$\frac{\sf 1}{\sf 5}$};
\node at (2/5, -.035) [scale=1, color=black] {$\frac{\sf 2}{\sf 5}$};
\node at (3/5, -.035) [scale=1, color=black] {$\frac{\sf 3}{\sf 5}$};
\node at (4/5, -.035) [scale=1, color=black] {$\frac{\sf 4}{\sf 5}$};
\node at (1/6, -.035) [scale=1, color=black] {$\frac{\sf 1}{\sf 6}$};

%\node at (0/3, 2/3) [scale=.9, color=black] {\sf b};
%\node at (1/2, 0/2) [scale=.9, color=black] {\sf c};

\draw [thick, red] (-1/2,0) circle (1/2+.618); %circle
\draw [thick, blue] (3/2,0) circle (1/2+.618); %circle

\end{tikzpicture}

\caption{The indicated points of tangency between the circles
are of form of fractions involving the Fibonacci numbers:
}
\label{fig:01}
\end{figure}

We finish with an image showing magnification of the circle system over the interval $[0,1]$
(Figure \ref{fig:01}).

\vspace{-.1in}
$$
\begin{array}{rl}
\hbox{\sf right sequence:} \quad &
\dfrac{1}{2}, \quad
\dfrac{3}{5}, \quad
\dfrac{8}{13}, \quad
...
\quad
\longrightarrow \  \dfrac{1}{\varphi} = \dfrac{\sqrt{5}-1}{2}
\\[12pt]
\hbox{\sf left sequence:} \quad &
\dfrac{1}{2}, \quad
\dfrac{2}{5}, \quad
\dfrac{5}{13}, \quad
...
\quad
\longrightarrow \  \dfrac{1}{\varphi^2} = \dfrac{3-\sqrt{5}}{2}
\end{array}
$$

%Bibliography-------------------------------------------------------------------------------

%?? [Add circle symbols to Ford circles to derive the coordinates of the tangency points. ]

%%%%%%%%%%%%%%%%%%%%%%%%%%%%%%%%%%%%%%
%%%%%%%%%%%%%%%%%%%%%%%%%%%%%%%%%%%%%%%

\begin{thebibliography}{99}

\bibitem 
{Apostol}
Apostol, T. M., Modular Functions and Dirichlet Series in Number Theory, 2nd ed. New York: Springer-Verlag, 1997
(§5.5, p. 101).

\vspace{-.1in}

\bibitem%
{Ford}
Lester Ford, Fractions, AMM, {\bf 45}, No. 9 (1938), pp. 586-601.

\vspace{-.1in}


\bibitem{jk}
Jerzy Kocik, A theorem on circle configurations, 	arXiv:0706.0372. 

\vspace{-.1in}

\bibitem
{Vajda}
Steven Vajda, Fibonacci and Lucas Numbers, and the Golden Section, Dover, 1989.

\vspace{-.1in}

\bibitem
{Weisstein}
Weisstein, Eric W. "Ford Circle." From MathWorld--A Wolfram Web Resource. http://mathworld.wolfram.com/FordCircle.html


\end{thebibliography}
\end{document}